\documentstyle{amsart}

\theoremstyle{plain}

\newtheorem{thm}{Theorem}[subsection]
\newtheorem{lem}{Lemma}[subsection]
\newtheorem{cor}{Corollary}[subsection]
\newtheorem{prop}{Proposition}[subsection]

\theoremstyle{definition}

\theoremstyle{remark}

\newtheorem{rem}{Remark}[subsection]
\newtheorem{ack}{Acknowledgments}  
\newtheorem{step}{Step}
\newtheorem{case}{Case}

\numberwithin{equation}{section}

\newsymbol\twoheadrightarrow 1310

\newcommand{\twr}{\twoheadrightarrow}
\newcommand{\riso}{\stackrel{\sim}{\longrightarrow}}
\newcommand{\kr}{\operatorname{ker}}
\newcommand{\rank}{\operatorname{rk}}
\newcommand{\codim}{\operatorname{codim}}
\newcommand{\pic}{\operatorname{Pic}}
\newcommand{\spec}{\operatorname{Spec}}

\newcommand{\Hom}{\operatorname{Hom}}
\newcommand{\ext}{\operatorname{Ext}}
\newcommand{\morph}{\operatorname{om}}
\newcommand{\col}{\colon}
\newcommand{\sm}{\setminus}
\newcommand{\im}{{\imath}}
\newcommand{\jm}{{\jmath}}
\newcommand{\vp}{\varphi}
\newcommand{\D}{\Delta}
\newcommand{\G}{\Gamma}
\newcommand{\lr}{\longrightarrow}

\newcommand{\ce}{{\cal E}}
\newcommand{\co}{{\cal O}}
\newcommand{\cs}{{\cal S}}
\newcommand{\cg}{{\cal G}}
\newcommand{\cf}{{\cal F}}
\newcommand{\cu}{{\cal U}}
\newcommand{\cv}{{\cal V}}

\newcommand{\cl}{{\cal L}}
\newcommand{\ct}{{\cal T}}
\newcommand{\cw}{{\cal W}}
\newcommand{\ch}{{\cal H}}
\newcommand{\cx}{{\cal X}}
\newcommand{\bbC}{{\Bbb C}}

\newcommand{\bbN}{{\Bbb N}}
\newcommand{\bbP}{{\Bbb P}}
\newcommand{\bbQ}{{\Bbb Q}}
\newcommand{\bbZ}{{\Bbb Z}}
\newcommand{\bbR}{{\Bbb R}}
\newcommand{\bbH}{{\Bbb H}}
\newcommand{\Q}{{\Bbb Q}}
\newcommand{\balp}{\boldsymbol{\alpha}}

\newcommand{\comp}{{\scriptstyle{\circ}}}

\newcommand{\End}{{\cal E}nd}
\newcommand{\sHom}{{\cal H}om}

\title{Intermediate Jacobians and Hodge Structures
 of Moduli Spaces}
\author{Donu Arapura}
\address{Department of Mathematics\\
Purdue University\\
West Lafayette, IN 47907\\
U.S.A.}
\thanks{First author partially supported by the NSF}
\email{dvb@@math.purdue.edu}
\author{Pramathanath Sastry}
\address{The Mehta Research Institute \\
Chhatnag, Jhusi, Allahabad District \\
U.P., INDIA 221 506}
\email{pramath@@mri.ernet.in}
\date{\today}
\begin{document}
\maketitle

{\em  Note: It was brought to our attention, by H. Esnault, that the hyperplane
$H$ in our thm 6.1.1 needs to be general, otherwise it's false. The counterexample
is to take $U$ to be the smooth part of a cone over an elliptic curve, and
$H$ to be a hyperplane passing through the vertex.
The flaw stems from the equation $j_!{\mathbb R}i'_*= {\mathbb R}i_*j'_!$ in the penultimate 
paragraph of its proof. As pointed out to us by M. Bondarko, this equation can be justified for 
general $H$ by an easy modification of the proof of condition $b_H$ on pp 35-36 of 
A. Beilinson, "How to glue perverse sheaves", Springer Lect Notes 128, (1987).
Note that the rest of the paper is unaffected by this change,
 since the theorem had only been  applied in the generic case. }

\section{Introduction}\label{s:intro}
We work throughout over the complex numbers $\bbC$, i.e. all schemes are over
$\bbC$ and all maps of schemes are maps of $\bbC$-schemes. A curve, unless
otherwise stated, is a smooth complete curve. Points mean geometric points.
We will, as is usual in such situations, toggle between the algebraic
and analytic categories without warning.

For a curve $X$, $\cs\cu_X(n,\,L)$ will denote the moduli space of
{\it semi-stable} vector bundles of rank $n$ and determinant $L$. The
smooth open subvariety defining the {\it stable locus} will be denoted
$\cs\cu_X^s(n,\,L)$. We assume familiarity with the basic facts about
such a moduli space as laid out, for example in
\cite[pp.\,51--52,\,VI.A]{css-drez} (see also
Theorems 10, 17 and 18 of {\it loc.cit.}).

When $L$ is a line bundle of degree coprime to $n$,
the moduli spaces $\cs\cu_X(n,L)$ and $\cs\cu_X^s(n,L)$ coincide,
and are therefore smooth and projective. The cohomology groups
$H^i(\cs\cu_X(n,L),\Q)$
carry  pure Hodge structures which can, in principle,
be determined by  using a natural set of  generators
(Atiyah-Bott \cite{at-bot}) and relations (Jeffrey-Kirwan
\cite{jef-kir}) for the cohomology ring; we will say more about this
later. When the degree of $L$ is not coprime to $n$ and
$g>2$, the situation is
complicated by the fact that
 $\cs\cu_X(n,L)$ is singular and
$\cs\cu_X^s(n,L)$ nonprojective. Thus the cohomology groups of
these spaces  carry (a-priori) mixed Hodge structures, and it are
these structures that we wish to understand. Our  main results concerns
the situation in low degrees.

\begin{thm} \label{thm:main1} Let 
$\i(n,g) = 2(n-1)g-(n-1)(n^2+3n+1)-7$.
Let $X$ be a curve of genus $g\ge 2$. If $n\ge 4$ and $i < \i(n,g)$ are 
integers, then for any pair of line bundles $L, L'$ (not necessarily
with the same degree) on $X$, the
mixed Hodge structures  $H^i(\cs\cu_X^s(n,L),\Q)$
and  $H^i(\cs\cu_X^s(n,L'),\Q)$ are
{\em (noncanonically)} isomorphic and  are both  pure of weight $i$.
\end{thm}

This statement is a bit disingenuous, it is vacuous unless $g\ge 16$.
The explicit determination of these Hodge structures is
rather delicate. However general considerations show that these
Hodge structures are semisimple and it is not
difficult   to write
down  all the potential candidates for the simple summands.

\begin{cor} With the notation as above, for $i < \i(n,\,g)$ any simple summand of
$$
H^i(\cs\cu^s_X(n,L),\Q)
$$ 
is, up to Tate twisting, a direct summand of
 a tensor power of $H^1(X)$
\end{cor}

For third cohomology, a more refined analysis yields:

\begin{thm}\label{thm:main} Let $X$ be a curve of genus $g\ge 2$, $n\ge 2$ an integer
and $L$ a  line bundle  on $X$.  Let $\cs^s=\cs\cu_X^s(n,\,L)$.
\begin{enumerate}
\item[(a)] If $g > {3\over n-1} + {n^2 + 3n + 3\over 2}$
and $n\ge 4$, then
$H^3(\cs^s,\bbZ)$ is a pure Hodge structure of type
$\{(1,2),\,(2,1)\}$, and it carries a
natural polarization making the intermediate Jacobian
$$J^2(\cs^s) = \frac{H^3(\cs^s,\,\bbC)}
{F^2+H^3(\cs^s,\,\bbZ)} $$
into a principally polarized abelian variety. There is an
isomorphism of principally polarized abelian varieties $J(X)\simeq
J^2(\cs^s)$.
\item[(b)] If $\deg L$ is a multiple of $n$, then the conclusions of 
{\em (a)}
are true for $g\ge 3$, $n\ge 2$ except the case $g=3,n=2$.
\end{enumerate}
\end{thm}

The word ``natural'' above has the following meaning: an isomorphism
between any two $\cs^s$'s as above will induce an isomorphism on third
cohomology  which will respect the indicated polarizations.
As an immediate corollary, we obtain the following  Torelli theorem:

\begin{cor}\label{cor:main} Let $X$ and $X'$ be curves of genus $g\ge 3$,
$L$ and $L'$ line bundles of (possibly different degrees)
 on $X$ and $X'$ respectively.
\begin{enumerate}
\item[(a)] Assume that
$n\ge 4$ is an integer such that $g > {3\over n-1} + {n^2 + 3n + 1\over
 2}$.
 If
\begin{equation}\label{eqn:stable}
\cs\cu_X^s(n,\,L) \simeq \cs\cu_{X'}^s(n,\,L')
\end{equation}
or if
\begin{equation}\label{eqn:semistable}
\cs\cu_X(n,\,L) \simeq \cs\cu_{X'}(n,\,L')
\end{equation}
then
$$
X \simeq X' .
$$
 
\item[(b)] If $\deg L = \deg L'$ and the common value is a multiple of
$n$, then the conclusions of {\em (a)} are true for $n\ge 2$, except the 
case $g=3,n=2$.
\end{enumerate}
\end{cor}

\begin{pf} Since $\cs\cu_X^s(n,\,L)$ (resp. $\cs\cu_{X'}^s(n,\,L')$)
is the smooth locus of $\cs\cu_X(n,\,L)$ (resp. $\cs\cu_{X'}(n,\,L')$),
therefore it is enough to assume \eqref{eqn:stable} holds.
By assumption $J^2(\cs\cu^s_X(n,\,L))\simeq J^2(\cs\cu^s_{X'}(n,\,L'))$
as polarized abelian varieties.
Therefore $J(X)\simeq J(X')$,
and the corollary follows from the usual Torelli theorem.
\end{pf}

When $(n,\,\deg\,L)=1$ (the ``coprime case"),
the second theorem (and its corollary with $\deg\, L = \deg\, L'$)
has been proven by Narasimhan and Ramanan \cite{N-R}, Tyurin
\cite{Ty} (both in the range $n\ge 2$ and $g\ge 2$, except when
$g=2, n=3$) and the special case of their results, when  $n=2$,
by Mumford and Newstead \cite{M-N}. In the non-coprime case, 
Kouvidakis and Pantev \cite{KP} had proved a Torelli theorem
for $\cs\cu_X(n,\, L)$, i.e. the  corollary 
under the assumption 
\eqref{eqn:semistable}, with better bounds.
In fact the full corollary can be deduced from this case.
However the present line of reasoning is extremely natural,
and is of a rather different character from that of Kouvidakis and 
Pantev. 
In particular,
Theorem\, \ref{thm:main} will not follow from their techniques.
In the special case where $n=2$ and $L= \co_X$, Balaji \cite{balaji}
has shown a similar Torelli type theorem
 for Seshadri's canonical desingularization
$N \to \cs\cu_X(2, \co_X)$ (see  \cite{css-desing}) in the range $g >3$.
\footnote{Balaji states the result for $g \ge 3$, but his proof seems to
work only for $g > 3$. (See Remark\,\ref{rmk:balaji}).}

Our strategy in the proof of both theorems is to use a Hecke
correspondence to relate the cohomology of $\cs\cu_X^s(n,L)$
with that of another moduli space $\cs\cu_X(n,L'')$ where the the
degree of $L''$ is coprime to $n$. When $n > 2$ the
maps defining the  Hecke correspondence are only  rational.
 And this 
necessitates some rather long calculations to bound the codimensions of the 
indeterminacy loci. Once the basic geometric properties of the correspondence 
are established,  the first theorem follows from some standard arguments in 
Hodge theory. For the second theorem, we need to make the isomorphism on
third cohomology canonical, and to moreover
 impose an intrinsic polarization on the Hodge structure
 $H^3(\cs\cu_X^s(n,\,L))$.

\section{The Main Ideas}\label{s:main-idea}
For the rest of the paper, we fix a curve $X$ of genus $g$, $n \in {\bbN}$,
$d \in {\bbZ}$  and a line bundle $L$ of degree
$d$ on $X$.  Let $\cs = \cs\cu_X(n,\,L)$ and $\cs^s =
\cs\cu^s_X(n,\,L)$. 

The main theorems will be proved in the final section of this
paper.
The broad strategy of our proofs are as follows\,:

\begin{step}

\begin{case}\label{step:one}

Assume, that $d$ is not divisible by $n$. 
Since $\cs\cu_X(n,\,L)$ is canonically
isomorphic to $\cs\cu_X(n,\,L\otimes\xi^n)$ for every line bundle
$\xi$ on $X$, we may assume that $0 < d < n$.

 Fix a set $\chi =
\{x^1,\ldots,\,x^{d-1}\} \subset X$ of $d-1$ distinct points.
Let
$$
\cs_1 =
\cs\cu_X(n,\,L\otimes\co_X(-D))
$$
where $D$ is the divisor
$x^1+\ldots +x^{d-1}$.

Construct (in \S\ref{s:hecke})  a generalized Hecke
correspondence consisting of a pair of rational maps
\begin{equation}
\cs_1 \stackrel{\scriptstyle{\pi}}{\leftarrow} {\bbP}
\stackrel{\scriptstyle{\phi}}
{\dashrightarrow}\cs^s
\end{equation}
By construction, there will be an  open subset $U\subset {\bbP}$
such that both $\pi|_U$ and $\phi|_U$ will be fiber bundles
 with fiber isomorphic to $({\bbP}^{n-1})^{d-1}$.
Estimates on the codimensions of the complements of $U$ and
 its image in $\cs^s$ will be
given in \S\ref{s:codim}. These estimates,
together with some  generalities on cohomology
and Hodge theory to be established in \S\ref{s:polar},
 will imply that
for small $i$, there are noncanonical isomorphisms of mixed
Hodge structures
\begin{equation}\label{eqn:hiSS}
H^i(\cs_1,\Q) \cong H^i(\cs^s, \Q).
\end{equation}
Therefore this reduces the proof of Theorem \ref{thm:main1} to
the case where $L$ and $L'$ have degree $1$, and this will be
treated in \S\,\ref{s:hodgedeg1}.
Moreover, for sufficiently large $g$, we have isomorphisms 
modulo torsion  of (integral, pure) Hodge structures
\begin{equation}\label{eqn:tate}
H^1(X,\,\bbZ)(-1) \riso H^3(\cs_1,\,\bbZ) \riso H^3(\cs^s,\,\bbZ),
\end{equation}
where the first isomorphism is the slant product with a certain
universal class (\S\ref{s:hodgedeg1}),
 and the second  now depends canonically on $X,\, L$
and $\chi$ (\S\ref{s:mainthms}). ``$(-1)$'' above is the Tate twist.
An isomorphism {\it modulo torsion} of integral pure Hodge structures means that
the underlying map of the finitely generated abelian groups is an isomorphism
on the free parts. In particular, if as above
these Hodge structure have odd weights,
the resulting map of the corresponding intermediate Jacobians is an isomorphism.

\end{case}

\begin{case} If $d$ is divisible by $n$ then we may assume that $d=0$.
In this case, set $\chi = \{x\}$ for some point $x\in X$.
Setting $D = x$, we construct $$
\cs_1 =
\cs\cu_X(n,\,L\otimes\co_X(-D))
$$
as before. In \S\ref{s:hecke0}, we construct a Hecke correspondence
analogous to the one above.
 However, now the map $\phi$ is regular
and  the codimension estimates are  substantially
better. This allows us to establish \ref{eqn:tate} with 
a much better bound on the genus.
\end{case}

\end{step}
\begin{step} Assume  that $g$ is
chosen sufficiently large. Our first task is to find a (possibly
nonprincipal) polarization $\Theta(\cs^s)$ on  $J^2(\cs^s)$ or
equivalently on the Hodge structure $H^3(\cs^s)$ which varies
algebraically with $X$.
The basic tools for constructing this are given in
\S\ref{s:polar}. Let
$$
\psi_{X,L,\chi}\col H^1(X)(-1) \riso
H^3(\cs^s)
$$
be the isomorphism given above,
and
$$\phi_{X,L,\chi}:J(X)\to J^2(\cs^s)$$
the corresponding isomorphism of abelian varieties.
 The isomorphisms vary algebraically  with the
data $(X,\,L,\,\chi)$. One can pull $\Theta(\cs^s)$ back
to get a second polarization on $J(X)$ which varies algebraically
with $(X,\,L,\,\chi)$. If we can find a positive  integer $m$ 
(independent of $(X,\,L,\,\chi)$) so that
$\Theta(\cs^s) = m\Theta$ where $\Theta$ is the standard
polarization, then it will follow that, after replacing $\Theta(\cs^s)$ with
$\frac{1}{m}\Theta(\cs^s)$, that $\Theta(\cs^s)$ is a principal
polarization such that
$(J(X),\Theta) \cong (J^2(\cs^s),\Theta(\cs^s))$ as
required. Since everything varies well, we can assume that
 $X$ is a sufficiently general curve in moduli. In this case,
 one checks that
any polarization on $J(X)$ is a multiple of the $\Theta$.
The precise argument is given in \S\ref{s:mainthms}.
\end{step}

\section{The Biregular Hecke Correspondence }\label{s:hecke}

The results of  this section  will be used to treat the case where 
$d=\deg L$ is not divisible by $n$. As explained earlier, we
may assume that $0< d < n$. 
We will continue the notation from step \ref{step:one}
of the previous section. 
  The degree of $L\otimes\co_X(-D)$
is $1$, therefore $\cs_1$ is smooth and there exists a Poincar{\'e}
bundle $\cw$ on $X\times\cs_1$. Let $\cw_1,\ldots,\,\cw_{d-1}$ be the
$d-1$ vector bundles on $\cs_1$ obtained by restricting $\cw$ to
$\{x^1\}\times{\cs_1}=\cs_1,\ldots,\,\{x^{d-1}\}\times{\cs_1}=\cs_1$
respectively. Let $\bbP_k=\bbP(\cw_k)$, $k=1,\ldots,\,d-1$,
where we use the convention $\bbP(W_i) = {\bf Proj}(S^*(W_i^*))$.
Let
$\bbP$ $(=\bbP_{X,L,\chi})$ be the product
$\bbP_1\times_{\cs_1}\ldots\times_{\cs_1}\bbP_{d-1}$.

\subsection{The map $\phi\col \bbP \dashrightarrow \cs$.}\label{ss:map-f}
 We need some notation\,:
\begin{itemize}
\item $\pi \colon\,\bbP \to \cs_1$ is the natural projection\,;
\item For $1\le k\le d-1$, $\pi_k\col \bbP \to \bbP_k$ is
the natural projection;
\item $\im\col Z \hookrightarrow X$ is the reduced subscheme defined
by $\chi=\{x^1,\ldots,\,x^{d-1}\}$.
\item $\im_k\col Z_k \hookrightarrow X$, the reduced scheme defined by
$\{x_k\}$, $k= 1,\ldots,\,d-1$.
\item For any scheme $S$,
\begin{enumerate}
\item[(i)] $p_S\col X\times\,S \to S$ and $q_S\col X\times\,S \to X$ are
the natural projections;
\item[(ii)] $Z^S = q_S^{-1}(Z)$;
\item[(iii)] $Z_k^S = q_S^{-1}(Z_k)$, $k=1,\ldots,\,d-1$.
Note that $Z^S_k$ can be identified canonically with $S$.
\end{enumerate}
\end{itemize}

We will show --- in \ref{ss:u-exact} --- that there is an exact sequence
\begin{equation}\label{eqn:u-exact}
0 \lr (1\times\pi)^*\cw \lr \cv \lr \ct_0 \lr 0
\end{equation}
on $X\times\,\bbP$, with $\cv$ a vector bundle on $X\times\bbP$
and $\ct_0$ (the direct image of) a line bundle {\it on the closed subscheme $Z^\bbP$},
which is universal
in the following sense\,: If $\psi\col S \to \cs_1$ is a $\cs_1$-scheme
and we have an exact sequence
\begin{equation}\label{eqn:s-exact}
0 \lr (1\times\psi)^*\cw \lr \ce \lr \ct \lr 0
\end{equation}
on $X\times\,S$, with $\ce$ a vector bundle on $X\times\,S$
and $\ct$ a line bundle {\it on the closed subscheme} $Z^S$, then
there is a unique map of $\cs_1$-schemes
$$
g\col S \lr \bbP
$$
such that,
$$
(1\times\,g)^*\eqref{eqn:u-exact} \equiv  \eqref{eqn:s-exact}.
$$
The $\equiv$ sign above means that the two exact sequences
are isomorphic, and the left most isomorphism
$(1\times{g})^*\comp(1\times{\pi})^* \riso (1\times\psi)^*$ is the canonical
one.

Let $U_1 \subset \bbP $ be the maximal open subset such that
$\cv|_{X\times\{t\}}$ is stable for each $t\in U_1$.
 We shall see that this is nonempty,
thus the natural moduli map $U_1\to \cs^s$ determines a rational
map  $\phi:\bbP{\dashrightarrow} \cs^s$. The geometric properties
of the map $\phi|_{U_1}$ are not obvious, so in \ref{ss:bireg}
we will  shrink $U_1$ to an open subset $U$ with more manageable properties.

\subsection{The universal exact sequence}\label{ss:u-exact}

We begin by reminding the reader of some elementary facts from commutative
algebra. If $A$ is a ring (commutative, with $1$), $t\in A$ a non-zero
divisor, and $M$ an $A$-module, then each element $m_0\in M$ gives
rise to an equivalence class of extensions
\begin{equation}\label{eqn:exact-mod}
0 \lr M \lr E_{m_0} \lr A/tA \lr 0
\end{equation}
where $E_{m_0} = \left(A\bigoplus{M}\right)/A(t,\,m_0)$, and the
arrows are the obvious ones. Moreover, if $m_0 - m_1 \in tM$, say
$$
m_0 - m_1 = tm'
$$
then the extension given by $m_0$ is equivalent to that given by $m_1$.
In fact, one checks that
\begin{equation}\label{eqn:patching}
\begin{split}
E_{m_0} & \riso E_{m_1} \\
(a,\,m) & \mapsto (a,\,m-am')
\end{split}
\end{equation}
gives the desired equivalence of extensions. This is another way of
expressing the well known fact that each element of
$M/tM = \ext^1(A/t,\,M)$ gives rise to an extension.

One globalizes to get the following\,: Let $S$ be a scheme,
$T\overset{\im}{\hookrightarrow} S$ a closed immersion, $\cf$ a
quasi-coherent $\co_S$-module, $U$ an open neighbourhood of $T$ in $S$,
and $t\in \G(U,\,\co_S)$ an element which defines $T \hookrightarrow U$,
and which is a non-zero divisor for $\G(V,\,\co_S)$
for any open $V \subset U$. Then every global section $s$ of
$\im^*\cf = \cf\otimes\co_T$ gives rise to an equivalence class of
extensions
\begin{equation}\label{eqn:patched}
0 \lr \cf \lr \ce \lr \co_T \lr 0.
\end{equation}
Indeed, we are reduced immediately to the case $S = U$. We build up
exact sequences \eqref{eqn:exact-mod} on each affine open subset
$W \subset S$, by picking a lift $\tilde{s_W}\in\G(W,\,\cf)$ of
$s\,|\,W$. One patches together these exact sequences via
\eqref{eqn:patching}.

Now consider $\bbP = \bbP_1\times_{\cs_1}\ldots\times_{\cs_1}\bbP_{d-1}$.
For each
$k = 1,\ldots,\,d-1$, let $p_k\col\,\bbP_k \to \cs_1$ be the natural
projection. We have a universal exact sequence
$$
0 \lr \co(-1) \lr p_k^*\cw_k \lr B \lr 0
$$
whence a global section $s_k\in\G(\bbP_k,\,p_k^*\cw_k(1))$.
However, note that
$$
p_k^*\cw_k = (1\times p_k)^*\cw\,|\,Z_k^{\bbP_k}
$$
where we are identifying $Z_k^{\bbP_k}$ with $\bbP_k$.
By \eqref{eqn:patched} we get exact sequences
$$
0 \lr (1\times\pi)^*\cw \lr \cv_k \lr
\co_{Z_k^\bbP}\otimes{L_k} \lr 0
$$
where $L_k$ is the line bundle obtained by pulling up $\co_{\bbP_k}(-1)$.
It is not hard to see
that $\cv_k$ is a family of vector bundles
parameterized by $\bbP$. Gluing these sequences together --- the $k$-th
and the $l$-th agree outside $Z^{\bbP}_k$ and $Z^{\bbP}_l$ --- we
obtain \eqref{eqn:u-exact}.

Now suppose we have a $\cs_1$-scheme $\psi\col\,S \to \cs_1$ and the
exact sequence \eqref{eqn:s-exact}
$$
0 \lr (1\times\psi)^*\cw \lr \ce \lr \ct \lr 0
$$
on $X\times{S}$. Restricting \eqref{eqn:s-exact} to $Z^S_k$ ($1\le k \le d-1$)
one checks that the kernel of $(1\times\psi)^*\cw\,|\,Z^S_k \to \ce\,|Z_k^S$
is a line bundle $\cl_k$. Identifying $Z^S_k$ with $S$, we see that
$(1\times\,\psi)^*\cw\,|\,Z^S_k = \psi^*\cw_k$. Thus $\cl_k$ is a line
sub-bundle of $\psi^*\cw_k$. By the universal property of $\bbP_k$, we
see that we have a unique map of $\cs_1$-schemes
$$
g_k\col\,S \lr \bbP_k
$$
such that $\co(-1)$ gets pulled back to $\cl_k$. The various $g_k$
give us a map
$$
g\col S \lr \bbP
$$
One checks that $g$ has the required universal property. The uniqueness
of $g$ follows from the uniqueness of each $g_k$.

\subsection {The biregular Hecke correspondence}\label{ss:bireg}

As explained earlier, we will need to consider an open subset
$U\subset U_1$, with good geometric properties. Essentially $U$
will be the preimage of a subset $U'\subset \cs^s$ parameterizing
those stable vector bundles $E$ for which the kernel of $p:E\to \co_Z$
is stable for every surjection. We again remind the reader, that if
$\deg L$ is a multiple of $n$, one can do better, as will be seen
in section \ref{s:hecke0}.

To construct $U$ and $U'$ rigorously, it would convenient to have
some sort of universal
family of vector bundles on $X\times \cs^s$. Unfortunately, such a
family will not exist when $(n,d) \not= 1$. To get around this, we
can go back to the construction (\cite{css-drez}): $\cs^s$ is
obtained as a good quotient of an open subscheme $R$ of a Quot scheme
by a reductive group $G$ (which is in fact a projective general linear
group). There is a vector bundle ${\Bbb E}'$ on $X \times R$
whose restriction to $X\times \{r\}$ is the bundle represented
by the image of $r$ in $\cs^s$. For $x\in X$, let $E_x$ be the restriction
of ${\Bbb E}'$ to $R \cong \{x\}\times R$. There is a natural $G$ action on
${\Bbb E}'$ which induces one on $E_x$ for every $x\in X$ and hence on
every $\bbP(E^*_x)$. This action is locally a product of the action on
$R$ with a trivial action along the fibers, therefore $\bbP(E_x^*)/G$
is a projective space bundle (or more accurately a Brauer-Severi
scheme) over $\cs^s$. Let $\pi':\bbP'\to R$ be the fiber product  of
$\bbP(E_{x^i}^*)\to R$ for $i=1,\ldots d-1$.
Then $\bbP'/G \to \cs^s$ is a $(\bbP^{n-1})^{d-1}$-bundle i.e.
a smooth morphism with fibers isomorphic to  $(\bbP^{n-1})^{d-1}$.
Let $L_i$ be the pullback of $\co_{\bbP(E_{x^i}^*)}(1)$ to $\bbP'$. 
There is a canonical map
$$
\kappa:{(1\times\pi')}^*{\Bbb E}' \to \bigoplus_k \,i_{k*}L_k
$$
where $i_k\colon\,\bbP' = \{x^k\}\times\bbP' \to X\times\bbP'$ is the natural
closed immersion.

Let $\cu''$ be the maximal open subset of $\bbP'$ such that 
$\kr(\kappa)|_{X\times\{t\}}$ is stable for each $t\in \cu''$.
Set $\cu' = R\sm\pi'(\bbP'\sm\cu'')$ and $\cu = {\pi'}^{-1}\cu'\subset
\cu''$. Both $\cu$ and $\cu'$ are invariant under $G$, and we let
$U$ and $U'\subset \cs^s$ be  the quotients. The map $U\to U'$ is
again a $(\bbP^{n-1})^{d-1}$-bundle. The exact sequence
$$
0\to \kr(\kappa)\to {(1\times\pi')}^*{\Bbb E}' \to \bigoplus_k \,i_{k*}L_k\to 0 
$$
yields, via the universal property of $\bbP$ a morphism $\cu\to \bbP$ which 
factors through $U$. The map $U\to \bbP$ is  an open immersion to a subset of $U_1$.
Thus we get a diagram
\begin{equation}\label{eqn:bireg-hecke}
\cs_1 \leftarrow U \stackrel{f}{\to} \cs^s
\end{equation}
which will be referred to as the {\it biregular Hecke correspondence}.

\begin{rem}\label{rem:fib-prod} With the above notation, note that ${\bbP}'/G \to
\cs^s$ is the fibre product of ${\bbP(E_{x^i}^*)/G}\to \cs^s$ ($i=1,\ldots,d-1)$.
It follows that $f$ itself may be regarded as the fibre product of $\bbP^{n-1}$
bundles over $U'$.
\end{rem}

\section{Codimension Estimates}\label{s:codim}

We will continue the notation from the previous section.
Our goal is to  establish the basic  estimates on the codimensions
of the complements of $U$ and $U'$. In particular, these
estimates will show that $U$ and $U'$ are nonempty.

\subsection { Special Subvarieties}

Let $X$ be a smooth curve of genus $g \ge 2$. Fix integers $n > m
\ge 1$, a rational number $\mu_0$ and a
line bundle $M$.
Let $T_{m,\mu_0}(n,M)$ be the subset of $\cs\cu^s(n,M)$
whose points  correspond to bundles $E$ for which there
exists a subbundle $F \subset E$ of rank $m$ and slope
$\mu( F) = \mu_0$.

\begin{lem} $T_{m,\mu_0}(n,M)$ is Zariski closed,
and can therefore be regarded as a reduced subscheme.
There exists a scheme $\Sigma$ such that $X\times \Sigma$
carries a rank $n$ vector bundle ${\Bbb E}$ with a 
rank $m$ subbundle ${\Bbb F}$. For each $\sigma\in \Sigma$,
the restriction of ${\Bbb E}$ to $X\times \{\sigma\}$ is
stable of degree $d$, while the restriction of ${\Bbb F}$ has
slope $\mu_0$.
The canonical map $\Sigma\to\cs\cu^s(n,M)$
has $T_{m,\mu_0}$ as its image.
\end{lem}

\begin{pf}
The argument is similar to that used above.
  $S=\cs\cu^s(n,M)$ can be realized as a good  quotient
of a subscheme $R$ of a Quot scheme by a reductive group $G$.
$X\times R$ carries a vector bundle ${\Bbb E}'$ whose restriction to
$X\times \{r\}$ is the bundle represented by the image of $r$ in $S$.
${\Bbb E}'$ extends to a coherent sheaf (denoted by the same symbol)
on the closure $Y = X\times \bar R$.
Let $Q$ be the closed subscheme
of  the relative Quot scheme 
$Quot_{X\times {\bar R}/ {\bar R}}({\Bbb E}')$ 
 parameterizing quotients which restrict to vector
bundles of degree $\deg M-\mu_0m$ and rank $n-m$ on each $X\times \{r\}$.
The the intersection $\tilde S$ of $R$ with the image of the
projection $ Q \to {\bar R}$ is closed
and $G$-invariant. Therefore the image of $\tilde S$ in
$S$, which is  $T_{m,\mu_0}$, is closed.

Set $\Sigma = R\times_{\bar R} Q$ and ${\Bbb E}$ to the pullback
of ${\Bbb E}'$ to $X\times \Sigma= (X\times R)\times_{\bar R} Q $.
Then it is easily seen
that these have the required properties.
\end{pf}

\subsection{ Deformations}

Let $D = \spec{\Bbb C}[\epsilon]/(\epsilon^2)$, and let $\ce$ be a
vector bundle on $X\times D$. $X$ can be  identified with a closed
subscheme of $X\times D$ with ideal sheaf ${\cal I} = \epsilon
O_{X\times D}$.
Let $E = \ce\otimes \co_{X} \cong \ce\otimes {\cal I}$. Then there is
an exact sequence
$$0 \to E \to \ce \to E \to 0$$
which is classified by an element of $Ext^1(E, E)\cong H^1(\End
(E))$. Furthermore
$$0 \to \det(E) \to \det(\ce) \to \det(E) \to 0$$
is classified by the trace of the above class. 
If $\End_0(E)$ denotes the traceless part of $\End(E)$, then
$H^1(\End_0(E))$ classifies  deformations of $E$
which induce trivial deformations of  $\det(E)$.

If $M= \det E$, then elements of $H^1(X, \End_0(E))$
give rise to maps from $D$ to $\cs\cu(n,M)$ which send the
closed point to the class of $E$. This map is  well
known to yield an isomorphism between   
$H^1(X, \End_0 (E))$ and the tangent space to $\cs\cu(n,M)$ at $[E]$.

Let $E$ be a vector bundle
corresponding to  general point $[E]$
of a component of $T_{m,\mu_0}$,
and let $v$ be a tangent vector to $T_{m,\mu_0}$ based
at $[E]$. This vector can be lifted to tangent vector
 $\tilde v$ to  $\Sigma$ at a point $\sigma$ lying over $[E]$.
Then $E\cong {\Bbb E}|_{X\times \{\sigma\}}$, and let
$F$ be the subbundle corresponding to 
${\Bbb F}|_{X\times \{\sigma\}} $. Set $G = E/F$,
$$K = \kr[\End (E) \to \sHom (F, G) ]$$
and
$$K_0 = K \cap \End_0(E) = \kr [\End_0 (E) \to \sHom (F, G) ].$$
 Note that $K = K_0\oplus (id)\cdot{\co_X}$.

Let $\ce$ be the first order deformation of $E$
corresponding to the tangent vector $v$.
Explicitly, if $D\to \Sigma$ is the map corresponding
to $\tilde v$, then $\ce = {\Bbb E}|_D$.
Setting ${\cal F} = {\Bbb F}|_D$ gives a deformation
of $F$ which fits into a diagram:
$$
\begin{array}{ccccccccc}
 0&  \to&      E& \to& \ce&        \to& E&         \to& 0\\
   &     &\uparrow&   & \uparrow &    & \uparrow&    &   \\
  0&  \to&      F& \to& {\cal F}&  \to& F&         \to& 0 .
\end{array}
$$

The images of the classes of ${\cal F}$ and $\ce$ in $Ext^1(F,E)$
agree up to sign. Therefore the class of $\ce$ lies in the kernel of the
map to $Ext^1(F,G)$ which is the image of $H^1(X, K_0)$.
Therefore:

\begin{lem} The image of 
$H^1(X,K_0)$ in $H^1(X,\End_0 (E))$
contains the tangent space to $T_{m,\mu_0}$ at $E$.
\end{lem}

A simple diagram chase shows that there is an exact sequence
$$0\to \sHom (E, F) \to K\to \End (G)\to 0$$
Therefore

\begin{cor}\label{cor:dimTbound}
  $\dim T_{m,\mu_0} \le h^1(\End (F)) +
h^1(\sHom(G,F)) + h^1(\End(G)) - g$.
\end{cor}

\subsection{Preliminary Codimension  Estimates}

 Let us say that $T_{m,\mu_0}$ is admissible
if  for any general point 
$[E]$ of a component of $T_{m,\mu_0}$ of largest dimension,
there exists a semistable subbundle
 $F\subset E$ of rank $m$ and slope $\mu_0$.  
We will estimate codimension of an admissible
$T_{m,\mu_0}(n,M)$ from below as a 
function of  four quantities $g\ge 2$,$n\ge 2$, $n>m\ge 1$, and
$\mu_0<\mu= n/deg\,M$. The imposition of admissibility 
simplifies the calculations, and presents no real loss of
generality.
 Fixing $E$ and $F$ as above,
let  $G = E/F$ and let
$0=G_0\subset G_1 \subset \ldots \subset G_r = G$ be
the Harder-Narasimhan filtration. Set $n_i = \rank(G_i/G_{i-1})$ and $n_0 = m$.
We have
$$
\mu <  \mu\left(G_{r}/G_{r-1}\right) < \cdots < \mu(G_1) < 
{n_0\over n_1}(\mu-\mu_0) + \mu
$$
where the last inequality follows from the bound on the slope of the
preimage of $G_1$ in $E$.

In the computations below, we will make use of the additivity of ``$\deg$", and
``$\rank$" and

\begin{lem} If $A$ and $B$ are locally free then $\mu(A\otimes B)
= \mu(A) + \mu(B)$ and $\mu(\sHom (A,B)) = \mu(B) - \mu(A)$.
\end{lem}

\begin {lem} If $V$ is semi-stable, then
 $h^0(V) \le \deg(V) + \rank(V)$ provided that the right side is
 nonnegative.
\end{lem}

\begin{pf}
This is obvious if $\deg V < 0$, since $h^0(V)=0$.
For $\deg V \ge 0$, we can assume that the lemma
holds for $V(-p)$ by induction. Therefore
$$h^0(V) \le h^0(V(-p)) + h^0(V/V(-p)) \le (\deg V - r + r) +r.$$

\end{pf}

\begin{cor} If $V$ is a semistable vector bundle then 
$h^1(\End (V)) \le \rank(V)^2g$.
\end{cor} 

\begin{pf}
As $\End (V)$ is semistable, we have
$h^0(\End (V)) \le \rank(V)^2$, and the result
now follows from Riemann-Roch. 
\end{pf}

Heuristically, $\dim T_{m,\mu_0}$ should be given by
the number of moduli for $F$, $G_i$
plus  extensions. To make this more rigorous,
we will work infinitesimally, and appeal to 
corollary \ref{cor:dimTbound}.
Each of the terms of the corollary
can be estimated in turn. For the first term, we have

\begin{equation}\label{eq:cod1}
h^1(\End (F)) \le n_0^2 g
\end{equation}

$Hom(G,F) =0$ by the numerical conditions, therefore
$$h^1(\sHom (G,F)) = -\chi(G^*\otimes F) $$
So  Riemann-Roch yields

\begin{equation}\label{eq:cod2}
 h^1(\sHom (G,F)) = n_0(n_1+\ldots n_r)(g-1) + nn_0(\mu - \mu_0).
\end{equation}
Note $\deg\,\sHom (G,F) = \deg\, \sHom (E, F) = -nn_0(\mu - \mu_0)$.

The last term $h^1(\End (G))$ remains. If $G=G_1$ is semi-stable, then
a bound is given  as above. To analyze the general case, we
need:

\begin {lem} If $V$ has a filtration such that the associated
graded sheaves are semi-stable bundle with slope at least $-1$, then
$h^1(V) \le g\cdot \rank(V)$.
\end{lem}

\begin{pf} By subadditivity of $h^1$, it is enough to assume that
$V$ is semi-stable with slope at least $-1$, in which case the result
follows from the previous lemma together with Riemann-Roch.
\end{pf}

\begin{cor} Let $W$ be a semi-stable bundle, and $V$ a vector
bundle with a  filtration such that the associated
graded sheaves are semi-stable bundle with slope at least $\mu(W)-1$, then
$$
h^1(\sHom (W, V)) \le \rank(V)\rank(W)g.
$$
\end{cor}

\begin{lem} We have the inequality
\begin{equation*}
h^1(\End (G_k)) \le (n_1 + \ldots n_k)^2g +
(\sum_{1\le i < j \le k} n_i n_j)\left( \mu(G_1)-\mu(G)-1\right).
\end{equation*}
\end{lem}

\begin{pf}
From
$$
0\to G_{k-1} \to G_{k} \to G_{k}/G_{k-1} \to 0,
$$
we obtain
\begin{align*}
h^1(\End (G_k)) \le  & h^1(\End (G_{k-1})) + h^1(\End (G_{k}/G_{k-1})) \\
 & + h^1(\sHom (G_{k-1}, G_{k}/G_{k-1})) + h^1(\sHom (G_{k}/G_{k-1}, G_{k-1})).
\end{align*}
It suffices to find upper bounds for each of the terms on the right,
and then sum them.
The first term can be controlled by induction, and the second and third 
terms by the previous  corollaries.
For the last term, an estimate can be obtained by  combining
$\Hom (G_{k}/G_{k-1}, G_{k-1}) = 0$, with
\begin{align*}
-\deg \sHom (G_{k}/G_{k-1}, G_{k-1}) &  = (n_1+ \ldots
n_{k-1})n_k\left[\mu(G_k/G_{k-1}) - \mu(G_{k-1})\right] \\
 &  = \sum_{i=1}^{k-1} n_i n_k\left[\mu(G_k/G_{k-1}) -\mu(G_i/G_{i-1})\right]\\
 & < \sum_{i=1}^{k-1} n_i n_k\left[\mu(G_1) -\mu\right]
\end{align*}
and the Riemann-Roch theorem.
\end{pf}

Thus
\begin{eqnarray}\label{eq:cod3}
h^1(\End (G_r)) - g &\le &  \left[(n_1 + \ldots n_r)^2-1\right]g\nonumber \\
&&    + (\sum_{1\le i < j \le r} n_i n_j)
\left(\mu(G_1)- \mu -1\right)
\end{eqnarray}

Subtracting equations \eqref{eq:cod1}, \eqref{eq:cod2} and \eqref{eq:cod3}
 from $\dim\, \cs\cu(n,M)$, simplifying, and
replacing $n_0$ by $m$,  yields
\begin{align*}
 \codim T_{m,\mu_0} > & m(n-m)(g+1) - (\sum_{1\le i < j \le r} n_i n_j)
\left(\mu(G_1)-\mu(G)-1\right) \\
 &  - nm(\mu-\mu_0) -(n^2-1) 
\end{align*}
To proceed further, we use
 $\sum_{1\le i < j \le r} n_i n_j \le {1\over 2}(n-m)^2$ and
$$
\mu(G_1) - \mu 
\le  {m \over n_1} (\mu - \mu_0) \le m(\mu-\mu_0)
$$
Putting these together yields

\begin{prop}\label{prop:codimT}
 If  $T_{m,\mu_0}$ is admissible, we have
$$\codim T_{m,\mu_0} >  m(n-m)(g+1)
- {1\over 2}m(n-m)^2(\mu - \mu_0)  
  -nm(\mu - \mu_0)- (n^2 -1).
$$
\end{prop}

\subsection {Final Codimension Estimates }

Recall that in \S\,\ref{s:hecke}  we had a diagram
$$
\cs_1\stackrel{\pi}{\leftarrow}
\bbP \supset U_1\supset U \stackrel{f}{\to} U'\subset \cs^s
$$
where $\bbP$ could be viewed as the parameter space for extensions
$$0 \to E \to E' \to \co_Z \to 0$$
with $E\in \cs_1$.
An extension $E'$ lies in $\bbP\sm U_1$ if and only if
 there exists a subbundle $F' \subset
 E'$ with $\mu(F') \ge \mu(E')= \mu(E) + (d-1)/n$.
$F'$ can be assumed to be semistable, since otherwise it can
be replaced by  the first step in the Harder-Narasimhan filtration.
 Let $F= F'\cap E$, then
$$
\mu(F) \ge  \mu(F') - {(d-1)\over m} \ge \mu(E) - 
(d-1)\left({1\over m}-{1\over n}\right)  \ge \mu(E) - {n-1\over m}
$$
where $m = \rank(F)$.

Similarly, if an extension  $E'$ lies in in $U_1\sm U$ then there exist a
coherent subsheaf $F_1' \subset E'$ which violates stability of the
kernel of some map $E'\to \co_Z$. After replacing $F_1'$ by 
the maximal semistable subbundle of its saturation, and
setting $F_1 = F_1' \cap E$, we obtain
$$
\mu(F_1) \ge \mu(F_1') - {(d-1)\over m'} \ge \mu(E) - {(d-1) \over m'}
\ge \mu(E) - {n-1\over m'}
$$
where   $m'= \rank(F_1)$.

Therefore $\pi(\bbP\sm U)$ is contained in the union of
 admissible $T_{m,\mu_0}(n,\,L\otimes\co_X(-D))$, as $m$ varies between $1$ and $n-1$,
and $\mu-\mu_0$ between zero and ${(n-1)/m}$.
Clearly $\codim(\bbP\sm U) \ge \codim(\pi(\bbP\sm U))$.
Therefore a term by term estimate of
the bound in Proposition\,\ref{prop:codimT}, using
$ (n-m) \le n-1$ and  $m(n-m)\ge n-1$ (for $1\le m\le n-1$),
 yields:

\begin{thm} $ \codim{(\bbP\setminus U)} > (n-1)\left[g - {1\over 2}(n^2 + n)\right] $
\end{thm}

\begin{cor}\label{cor:keyest1}
If 
$$
g > {k\over n-1} + {n^2 +n\over 2} 
$$ 
then,
$$
\codim(\bbP\setminus U) > k .
$$
\end{cor}

The estimate on the codimension of $U'$ is obtained by a similar
argument. If $E' \in \cs\sm U'$, then there exists a surjection
$E' \to \co_Z$ and a subbundle $F\subset \kr[E'\to \co_Z]$ violates the
stability of the kernel. Let $m = \rank(F)$ and
 $F'$ be  
the maximal semistable sub-bundle of the saturation of $F$ in $E'$, then
\begin{align*}
\mu(F) & \ge \mu(F') - {(d-1)\over m} \\
 & \ge \mu(E) - {(d-1) \over m} \\
 & = \mu(E') - (d-1)\left({1\over n} + {1\over m}\right) \\
 & \ge \mu(E') - (n-1)\left({1\over n} + {1\over m}\right) .
\end{align*}
Thus the complement of $U'$ lies in the union of
admissible 
$T_{m,\mu_0}(n,\,L)$ as $m$ varies between $1$ and $n-1$,
and $\mu-\mu_0$ varies between zero and 
$(n-1)({1\over n} + {1\over m})$. 
An elementary analysis shows that
$$
\left(1 + {m\over n}\right)\left[{1\over 2}(n-m)^2 + n\right] \le 
\left(1+ {1\over n}\right)\left[{1\over 2}(n-1)^2 + n\right]
$$
for $1\le m \le n-1$ and $n \ge 4$.
Combining  this
and the trivial estimate $m(n-m)\ge n-1$
 with Proposition\,\ref{prop:codimT} yields:

\begin{thm}
If $n \ge 4$ then
$$
\codim(\cs^s\setminus U) > (n-1)\left[g - {n^2 +3n+1\over 2} -
{3\over n}\right]
> (n-1)\left[g-{n^2 +3n+1\over 2}\right] -3 
$$
\end{thm}

\begin{cor}\label{cor:keyest2}
 If $n \ge 4$ and $g > {k\over n-1} + {n^2 + 3n + 3\over 2} $,
 then $\codim(\cs^s\sm U') > k$.
\end{cor}

Of course, the above inequality for $g$ implies the inequality
in corollary \ref{cor:keyest1}.
The restriction on $n$ above is harmless, because  the
cases of $n=2,\, 3$ can be handled by the results of the
next section.

\section{Hecke When $\deg L = 0$}\label{s:hecke0}

The notations in this section are special to this section.
We now give a an improved Hecke correspondence when $\deg L\in n\bbZ$.
Clearly, without loss of generality, we may (and will) assume that 
$\deg L =0$.  In this case, instead of choosing $d-1$ points on $X$ (which 
clearly does not make sense), we choose one point $x\in X$, and let $Z$ be 
the reduced scheme supported on $\{x\}$. Let $D$ be 
the divisor given by $\{x\}$. Then $\deg{L\otimes\co_X(-D)}$ is $-1$, and if
$\cs_1 = \cs\cu_X(n,\,L\otimes\co_X(-D))$, then $\cs_1$ is smooth and
there exists a Poincar\'e bundle $\cw$ on $X\times\cs_1$. Let $\cw_1$
be the vector bundle on $\cs_1$ obtained by restricting $\cw$ to
$\{x\}\times \cs_1 = \cs_1$, and $\bbP = \bbP(\cw_1)$. Let $\pi:\bbP \to \cs_1$
be the natural projection. Then as before we have the universal exact sequence 
\eqref{eqn:u-exact} 
$$
0 \to (1\times\pi)^*\cw \to \cv \to \ct_0 \to 0,
$$
of coherent sheaves on $X\times\bbP$ with $\cv$ a vector bundle, and
$\ct_0$ a sheaf supported on $\{x\}\times\bbP = \bbP$ which is a line
bundle on $\bbP$. This means that $\bbP$ parameterizes exact sequences
$$
0 \to W \to V \to \co_Z \to 0
$$
of coherent sheaves on $X$ with $W\in\cs_1$ and $V$ a vector bundle
(necessarily of rank $n$ and determinant $L$).

There is a way of interpreting this universal property in terms of
quasi-parabolic bundles (see \cite[p.\,211--212,\,Definition\,1.5]{mehta-css}
for the definitions of quasi-parabolic and parabolic bundles).  
We introduce a quasi-parabolic datum on $X$ by attaching the flag type 
$(1,\,n-1)$ to the point $x$. From now onwards {\it quasi-parabolic structures
will be with respect to this datum and on vector bundles of rank $n$ and
determinant $L$}.  
One observes that for a vector bundle $V$ (of rank $n$ and determinant
$L$), a surjective map $V\twr \co_Z$ determines a unique quasi-parabolic
structure, and two such surjections give the same quasi-parabolic structure
if and only if they differ by a scalar multiple. 
The above mentioned universal property says that $\bbP$
is a (fine) moduli space for quasi-parabolic bundles.  More precisely, the
family of quasi-parabolic structures
$$
\cv \twr \ct_0
$$
parameterized by $\bbP$ is universal for families of quasi-parabolic bundles
$$
\ce \twr \ct
$$
parameterized by $S$, whose kernel is a family of semi-stable bundles.  The
points of $\bbP$ parameterize quasi-parabolic structures $V \twr \co_Z$
whose kernel is semi-stable. 

Let $\balp = (\alpha_1,\,\alpha_2)$, where $0 < \alpha_1 < \alpha_2 <1$,
and let $\D = \D_{\balp}$ be the parabolic datum which attaches 
weights $\alpha_1, \alpha_2$ to our quasi-parabolic datum.
We can choose $\alpha_1$ and $\alpha_2$ so small that
\begin{itemize}
\item a parabolic semi-stable bundle is parabolic stable\,;
\item if $V$ is stable, then every parabolic structure on $V$ is parabolic
stable\,;
\item the underlying vector bundle of a parabolic stable bundle is
semi-stable in the usual sense.
\end{itemize}
Showing the above involves some very elementary calculations. Call $\balp$
{\em small} if $\alpha_1$ and $\alpha_2$ satisfy the above properties and
denote the resulting (fine) moduli space of parabolic stable bundles 
$\cs\cu_X(n,\,L,\,\D)$. 

\begin{lem}\label{lem:ker-ss} If $\balp$ is small then for every parabolic
stable bundle $V\twr\co_Z$, the kernel $W$ is semi-stable.
\end{lem}

\begin{pf} Note that $\deg{W} = -1$ and hence the semi-stability of $W$ is
equivalent to its stability. Suppose $W$ is not stable. Then by the above
observation, there is a subbundle $E$ of $W$ such that $\mu(E) > \mu(W)$.
Let $\rank{E} = r$. Let $E'\subset V$ be the subbundle generated by $E$. 
Let $T$ be the torsion subsheaf of the cokernel of $E \to V$, and $t$ the
vector space dimension (over $k(x)=\bbC$) of $T$. We then have
$\deg{E'} = \deg E + t$, where $t\ge 0$. Thus
$$
-\dfrac{1}{n} < \mu(E) \le \mu(E') = \mu(E) + \dfrac{t}{r} \le \mu(V) = 0.
$$
In particular
\begin{align*}
-\dfrac{1}{n} & < \mu(E) \le -\dfrac{t}{r}, \\
\intertext{i.e.}
\dfrac{1}{n} & > \dfrac{t}{r},
\end{align*}
but this is not possible for $t > 0$ since $r < n$. So $t=0$ and hence $T=0$.
Thus $E = E'$. Let $d = \deg{E}$. Then
$$
-\dfrac{1}{n} < \dfrac{d}{r} \le 0.
$$
This is possible only if $d=0$. Since $E=E'$, one checks that the flag
induced on the fibre $E'_x$ by the flag $F_1V_x \supset F_2V_x$ is the
trivial flag $E'_x = F_1E'_x = F_2E'_x$. It follows that the parabolic
degree of $E'$ is $r\alpha_2$. Since $V$ is parabolically stable, this
means
$$
\alpha_2 < \dfrac{\alpha_1 + (n-1)\alpha_2}{n}.
$$
The right side is a non-trivial convex combination of $\alpha_1$ and
$\alpha_2$, and we also have $\alpha_1 < \alpha_2$. Thus we have
a contradiction. Therefore $W$ is stable (see also \cite{balaji-thesis}
for the same result when $n=2$).
\end{pf}

\begin{thm}\label{thm:par-moduli} $\bbP = \cs\cu_X(n,\,L,\,\D)$, and
$\cv \twr \ct_0$ is the universal family of parabolic bundles.
\end{thm}

\begin{pf}
Let $\bbP^{ss}\subset \bbP$ be the locus on which $\cv$ consists
of parabolic semi-stable (=parabolic stable) bundles.  One checks that
$\bbP^{ss}$ is an open subscheme of $\bbP$ (this involves two things\,:
(i) knowing that the scheme $\widetilde{R}$ of \cite[p.\,226]{mehta-css}
has a local universal property for parabolic bundles and (ii) knowing that
the scheme ${\widetilde{R}}^{ss}$ of {\it loc. cit. } is open). 

Clearly $\bbP^{ss}$ is non-empty --- in fact if $V$ is stable of rank $n$
and determinant $L$, then any parabolic structure on $V$ is parabolic
stable (see above).  We claim that $\bbP^{ss} \simeq \cs\cu_X(n,\,L,\,\D)$. 
To that end, let $S$ be a scheme, and
\begin{equation}\label{eqn:par-fly}
\ce \twr \ct
\end{equation}
a family of parabolic stable bundles parameterized by $S$. By 
Lemma\,\ref{lem:ker-ss}, the kernel
$\cw'$ of \eqref{eqn:par-fly} is a family of stable bundles of rank $n$
and determinant $L\otimes\co_X(-D)$.  Since $\cs_1$ is a fine moduli space,
we have a unique map $g\col S\to \cs_1$ and a line bundle $\xi$ on $S$
such that $(1\times{g})^*\cw = \cw'\otimes\,p_S^*\xi$.  By doctoring
\eqref{eqn:par-fly} we may assume that $\xi=\co_S$.  The universal
property of the exact sequence \eqref{eqn:u-exact} on $\bbP$ then
gives us a unique map
$$
g\col\,S \lr \bbP
$$
such that $(1\times{g})^*\eqref{eqn:u-exact}$ is equivalent to
$$
0 \lr \cw' \lr \ce \lr \ct \lr 0. 
$$
Clearly $g$ factors through $\bbP^{ss}$.  This proves that
$\bbP^{ss}$ is $\cs\cu_X(n,\,L,\,\D)$.  However, $\cs\cu_X(n,\,L,\,\D)$
is a projective variety (see \cite[pp.\,225--226,\,Theorem\,4.1]{mehta-css}),
whence we have 
$$
\bbP = \cs\cu_X(n,\,L,\,\D). 
$$
Clearly, $\cv\twr\ct_0$ is the universal family of parabolic bundles.
\end{pf}

Note that  the above proof gives $\bbP^{ss} = \bbP$, whence we have,

\begin{cor} Let 
$$
0 \to W \to V \to \co_Z \to 0
$$
be an exact sequence of coherent sheaves on $X$ with $W\in \cs_1$ and $V$
a vector bundle. Then $V$ is a semi-stable bundle.
\end{cor}

It follows that $\cv$ consists of 
(usual) semi-stable bundles (by our choice of $\balp$).  Since $\cs$
is a coarse moduli space, we get the map
\begin{equation}\label{eqn:phi}
\vp\col\,\bbP \lr \cs . 
\end{equation}

\begin{rem}\label{rmk:hecke}
Note that the parabolic structure $\D$ is something of a red herring. 
In fact $\cs\cu_X(n,\,L,\,\D)$ parameterizes quasi-parabolic structures
$V \twr \co_Z$, whose kernel is semi-stable
(cf.  \cite[p.\,238,\,Remark\,(5.4)]{mehta-css} where this point is made
for $n=2$). 
\end{rem}

\begin{rem}\label{rem:proj-bundle} 
Let $V$ be a stable bundle of rank $n$, with $\det{V}=L$, so that
(the isomorphism class of) $V$ lies in $\cs^s$.  Since any parabolic
structure on $V$ is parabolic stable (by our choice of $\balp$), therefore
we see that $f^{-1}(V)$ is canonically isomorphic to
$\bbP(V_{x}^*)$\footnote{One can
be more rigorous.  Identifying $Z^{\bbP} = \{x\}\times\bbP$ with ${\bbP}$ 
we see that restricting the universal exact sequence to $Z^{\bbP}$ gives us 
the quotient $\co_{\bbP}\otimes_{\bbC}V_{x} \twr \ct_0|Z^{\bbP}$.  Let $S$
be a scheme which has a quotient $\co_S\otimes_{\bbC}V_{x}\twr \cl$ on it
where $\cl$ is a line bundle. This quotient extends (uniquely) to a family
parabolic structures $q_S^*V \twr \ct$ (on $V$) parameterized by $S$. 
By the Lemma, the kernel is a family of stable bundles. The
universal property of the exact sequence \eqref{eqn:u-exact} gives us a map
$S \to \bbP$, and this map factors through $f^{-1}(V)$. }. Moreover, 
it is not hard to see that $\bbP^s := \pi^{-1}(\cs^s) \to \cs^s$ is smooth 
(examine the effect on the tangent space of each point on $\bbP^s$). 
\end{rem}
\subsection{Codimension estimates. }\label{ss:codim}  We wish to estimate 
$\codim{(\bbP\sm \bbP^s})$ as well as $\codim{(\cs\sm\cs^s)}$. The second admits
to exact answers (see Remark\,\ref{rmk:codim} below). Heuristically (one
can make this rigorous via the deformation theoretic techniques in
\S\,\ref{s:codim}), the method for obtaining the first estimate is as follows.   

Let $V \twr \co_Z$ be a parabolic bundle in $\bbP\sm \bbP^s$.  Then we have 
a filtration (see \cite[p.\,18,Th\'eor\`eme\,10]{css-drez})
$$
0=V_{p+1} \subset V_p \subset \ldots \subset V_0=V
$$
such that for $0\le i \le p$, $G_i=V_i/V_{i+1}$ is stable and $\mu(G_i)=\mu$. 
Moreover (the isomorphism class of) the vector bundle $\bigoplus{G_i}$ depends
only upon $V$ and not on the given filtration.  We wish to count the
number of moduli at $[V\overset{\theta}{\twr} \co_Z] \in \bbP\sm \bbP^s$. 
There are three sources\,:
\begin{enumerate}
\item[a)] The choice of $\bigoplus_{i=0}^pG_i$\,;
\item[b)] Extension data\,;
\item[c)] The choice of parabolic structure $V\overset{\theta}{\twr}\co_Z$,
for fixed semi-stable $V$. 
\end{enumerate}
The source c) is the easiest to calculate --- there is a codimension
one subspace at each parabolic vertex, contributing $n-1$ moduli. Let 
$n_i=\rank{G_i}$. The number of moduli arising from a) is evidently
$$
\sum_{i=0}^p(n_i^2-1)(g-1) + pg . 
$$
Indeed, the bundles $G_i$ have degree $n_i\mu$ and the product of their
determinants must be $L$.  They are otherwise unconstrained. For c),
again using techniques in \S\,\ref{s:codim}, one sees that the number
of moduli contributed by extensions is
\begin{equation*}
\begin{split}
\sum_{i=0}^p\left [h^1(G_i^*\otimes{V_{i+1}}) -1 \right]  & \le
\sum_{i=0}^p\left[p-i-n_i(n_{i+1}+\ldots n_p)(1-g) \right] - (p+1) \\
 & = \dfrac{p(p+1)}{2} - 
\sum_{i=0}^{p-1}n_i(n_{i+1}+\ldots +n_p)(1-g) - (p+1) \\
 & = \dfrac{(p+1)(p-2)}{2} - \sum_{i<j}n_in_j(1-g). 
\end{split}
\end{equation*}
This gives
\begin{equation*}
\begin{split}
\codim(\bbP\sm \bbP^s) & \ge \sum_{i<j}n_in_j(g-1) - \dfrac{(p-1)(p+2)}{2} \\
& = B \qquad (\text{say}). 
\end{split}
\end{equation*}
Now,
$\sum_{i<j}n_in_j \ge {p(p+1)}/{2}$, therefore
$$
B \ge \dfrac{p(p+1)}{2}(g-1) - \dfrac{(p+2)(p-1)}{2}. 
$$
It follows that $B\ge 3$ whenever $p\ge 2$ {\it and} $g\ge 3$.  If $p=1$ and
$n \ge 3$, then 
$$
B/(g-1) = \sum_{i<j}n_in_j \ge 2
$$
and one checks that $B\ge 3$ whenever $g\ge 3$.  

\begin{rem}\label{rmk:codim} We could use similar techniques
to estimate $\codim{(\cs\sm \cs^s)}$, but our task is made
easier by the exact answers in \cite[p.\,48,\,A]{css-drez}. For just
this remark, assume $d > n(2g-1)$, and let $a = (n,\,d)$.  Then
$a\ge 2$.  Let $n_0 = n/a$.  Then according to {\it loc. cit. },
\begin{equation*}
\codim{(\cs\sm \cs^s)}=
\begin{cases}
(n^2-1)(g-1) -
\dfrac{n^2}{2}(g-1) -2 +g & \text{if $a$ is even} \\
 & {} \\
(n^2-1)(g-1) - \dfrac{n^2 + n_0^2}{2}(g-1) -2 + g & \text{if $a$ is odd}. 
\end{cases}
\end{equation*}
It now follows that
$$
\codim{(\cs\sm \cs^s)} > 5
$$
if $n,\,g$ are in the range of Theorem\,\ref{thm:main}(b).
\end{rem}
%
\section {Hodge theory}\label{s:polar}

This section, which can be read independently of the rest of the
paper, contains some results from Hodge theory that
will be  needed to complete the proofs of the main theorems.

\subsection{Purity}

We refer to  \cite{deligne-hodge} for the definition and
basic properties of mixed Hodge structures.
Deligne's fundamental result is that the cohomology groups
of  schemes with coefficients in $\bbZ$
carry canonical mixed Hodge structures.
We will need certain purity results for these mixed Hodge structure
for low degree cohomology of smooth open varieties. These results
can be deduced by comparing ordinary cohomology to intersection
cohomology and appealing to the work Saito \cite{saito}.
However we will give more elementary arguments, using a version of
the Lefschetz hyperplane theorem.

The notation of this section will be independent of the others.
\begin{lem}\label{lem:codim} If $Y$ is a smooth  variety,
$Z$ a codimension $k$ closed subscheme, and $U=Y\sm Z$, then
$$
H^j(Y,\,\bbZ) \riso H^j(U,\,\bbZ)
$$
for $j < 2k-1$.
\end{lem}
\begin{pf}
We have to show that $H^j_Z(Y,\,\bbZ)$ vanishes for $j < 2k$. By Alexander
duality (see for e.g. \cite[p.\,381,\,Theorem\,4.7]{iverson}) we have
$$
H^j_Z(Y,\,\bbZ) \riso H_{2m-j}(Z,\,\bbZ),
$$
where $m=\dim{Y}$ and $H_*$ is Borel-Moore homology. Now use
\cite[p.\,406,\,3.1]{iverson} to conclude that the
right side vanishes if $j < 2k$ (note that $``\dim"$ in
{\it loc.cit} is dimension as an analytic space,
and in {\it op.cit.} it is dimension as a topological (real)
manifold).
 \end{pf}

\begin{rem}\label{rmk:balaji} In view of the above Lemma, it
seems that Balaji's proof of Torelli (for Seshadri's
desingularization of $\cs\cu_X(2,\,\co_X)$) does not work
for $g=3$, for in this case, the codimension of $\bbP\sm \bbP^s =2$. 
(See \cite[top of p.\,624]{balaji} and \cite[Remark\,9]{balaji-thesis}.)
\end{rem}

\begin{cor} $H^j(U,\,\bbZ)$ is pure of weight $j$ when $j< 2k-1$.
\end{cor}

We will need purity results even when compactification is not
smooth. To this end we prove the following version of the
Lefschetz theorem.

\begin{thm}\label{thm:lefschetz}
Let $Y$ be an $m$-dimensional projective variety.
Suppose that $U$ is a smooth Zariski open subset.
 If $H$ is a hyperplane section
of $Y$ such that $U\cap H$ is non-empty, then
$$
H^i(U,\,\bbQ) \to H^i(U\cap{H},\,\bbQ)
$$
is an isomorphism for $i < m-1$ and injective when $i = m-1$.
\end{thm}
\begin{pf} We need some results involving Verdier duality. The
standard references are \cite{borel} and \cite{iverson}. Let
$S$ be an analytic space and $p_S$ the map from $S$ to a point.
For $\cf\in D^b_{const}(S,\,\bbQ)$ (the derived category of
bounded complexes of $\bbQ_S$-sheaves whose cohomology sheaves
are $\bbQ_S$-constructible), set
$$
D_S(\cf) = \bbR\ch\morph_S(\cf,\,p_S^!\bbQ).
$$
We then have by Verdier duality
\begin{equation}\label{eqn:verdier}
\bbH^i(S,\,\cf) \overset{\sim}{\lr} \bbH^{-i}(S,\,D_S(\cf))^*.
\end{equation}
Here $\bbH^*$ denotes ``hypercohomology".

For an open immersion $h\col S' \hookrightarrow S$, one has canonical
isomorphisms
\begin{align}\label{eqn:*!}
\bbR{h}_*D_{S'}\cg & \overset{\sim}{\lr} D_S(h_!\cg) \\
\intertext{and}
\label{eqn:!*}
{\bbR}h_!D_{S'}\cg & \overset{\sim}{\lr} D_S(\bbR{h_*}\cg)
\end{align}
Here $\cg\in D^b_{const}(S',\,\bbQ)$.
The first isomorphism is easy (using Verdier duality for the
map $h$) and the second follows from the first and from the
fact that $D_{S'}$ is an involution. We have used (throughout)
the fact that $h_!$ is an exact functor.

If $S$ is smooth, we have
\begin{equation}\label{eqn:dim}
p_S^!\bbQ = \bbQ_S[2\dim{S}].
\end{equation}
In order to prove the theorem, let $V=U\sm H$ and $W=Y\sm H$. We then
have a cartesian square
$$
\begin{array}{ccc}
V & \stackrel{\scriptstyle{{\im}'}}{\lr} & U \\
\vcenter{%
\llap{$\scriptstyle{{\jm}'}$}}\Big\downarrow &
& \Big\downarrow\vcenter{%
\rlap{$\scriptstyle{\jm}$}} \\
W & \underset{\scriptstyle{\im}}{\lr} & Y
\end{array}
$$
where each arrow is the obvious open immersion.
We have, by \eqref{eqn:*!} and \eqref{eqn:!*}, the identity
\begin{equation}\label{eqn:Dji}
{\jm}_!\bbR{\im}'_*D_V\bbQ_V = D_Y(\bbR\jm_*{\im}'_!\bbQ_V).
\end{equation}
Consider the exact sequence of sheaves
$$
0 \lr {\im}'_!\bbQ_V \lr \bbQ_U \lr g_*\bbQ_{H\cap U} \lr 0
$$
where $g\col H\cap{U}\to U$ is the natural closed immersion. It suffices
to prove that $H^i(U,\,{\im}'_!\bbQ_V)=0$ for $i\le m-1$. Now,
$$
H^i(U,\,{\im}'_!\bbQ_V)=\bbH^i(Y,\,\bbR\jm_*{\im}'_!\bbQ_V).
$$
Using \eqref{eqn:verdier}, \eqref{eqn:Dji} and \eqref{eqn:dim}, the
above is dual to
\begin{align*}
\bbH^{-i}(Y,\,\jm_!\bbR{\im}'_*D_V\bbQ_V)
 & = \bbH^{2m-i}(Y,\,\jm_!\bbR{\im}'_*\bbQ_V) \\
\intertext{But $\jm_!\bbR{\im}'_* = \bbR\im_*{\jm}'_!$, and hence the above is}
 & = \bbH^{2m-i}(Y,\,\bbR{\im}_*({\jm}'_!\bbQ_V)) \\
 & = \bbH^{2m-i}(W,\,{\jm}'_!\bbQ_V) \\
 & = H^{2m-i}(W,\,{\jm}'_!\bbQ_V).
\end{align*}
Now, $W$ is an affine variety, and therefore, according to
M. Artin, its constructible cohomological
dimension is less than or equal to its dimension \cite{artin}.
Consequently, the above chain of equalities vanish whenever
$i<m$ (see also \cite{gor-mac}).
\end{pf}

We will use the notation of this theorem for the remainder of the
section. We immediately have:
\begin{cor}\label{cor:lefschetz} Let $e=\codim(Y\sm U)$.
For $i< e-1$, the Hodge structure
$H^i(U,\,\bbZ)$ is pure of weight $i$.
\end{cor}
\begin{pf} This is true if $U$ is projective. In general proceed
using Bertini's theorem, induction, Theorem\,\ref{thm:lefschetz} and
the fact that submixed Hodge structures of pure Hodge structures are
pure \cite{deligne-hodge}.
\end{pf}

\subsection{Polarizations.} Let $Y, U, e$, etc. be as in Theorem\,\ref{thm:lefschetz}.
Let $i\in\bbN$ and $\cl$ a line bundle on $Y$ be such that
\begin{enumerate}
\item[(a)] $H^j(U,\,\bbQ)=0$ for $j = i-2,\,i-4,\ldots $ (note that this forces
$i$ to be {\it odd})\,;
\item[(b)] $i < e-1$\,;
\item[(c)] $\cl$ is very ample.
\end{enumerate}

Let $M$ be the intersection of $k=m-e+1$  hyperplanes in general position.
Then $M$ is a smooth variety contained in $U$. Let
$$
l\col H^i(U) \lr H^{2m-i}_c(U)
$$
be the composite of
\begin{equation*}
\begin{split}
H^i(U) & \lr H^i(M) \\
 & \lr H^{2m-2k-i}(M) \\
 & \lr H^{2m-i}_c(U)
\end{split}
\end{equation*}
where the first map is restriction, the second is
``cupping with $c_1(\cl)^{m-k-i}$" and the third is the Poincar{\'e} dual
to restriction. The map $l$ is also described as
$$
x \mapsto x\cup c_1(\cl)^{m-k-i}\cup [M].
$$
One then has (easily)
\begin{lem} If $M'$ is another $k$-fold intersection of general hyperplanes,
then $[M'] = [M]$. Therefore $l$ depends only on $\cl$.
\end{lem}
\begin{prop}\label{prop:polarization} The pairing
$$
<x,\,y> = \int_Ul(x)\cup y
$$
on $H^i(U,\,\bbC)$
gives a polarization on the pure Hodge structure $H^i(U,\,\bbZ)$.
This  makes the associated
complex torus $J^p(U)$ (where $i = 2p-1$) into an abelian variety
when $H^i(U)$ is of type $\{(p,p-1),(p-1,p)\}$.
\end{prop}
\begin{pf} By Theorem\,\ref{thm:lefschetz}, we have an isomorphism
$$
r\col H^i(U) \lr H^i(M).
$$
The latter Hodge structure carries a polarization given by
$$
<\alpha,\,\beta> = \int_Mc_1(\cl)^{m-k-i}\cup\alpha\cup\beta.
$$
The conditions on $i$ and the Hodge-Riemann bilinear relations
on the primitive part of $H^i(M,\,\bbC)$, assure us that the above is indeed
a polarization (see \cite{grif-period} or \cite[Chap.\,V,\,\S6]{wells}).
Our conditions on $i$ imply that the primitive part of $H^i(M)$ is everything
so that the $p$-th intermediate Jacobian of $M$ is an abelian variety. 
The pairing on $H^i(M)$ translates to a polarization on $H^i(U)$ given by
$$
<x,\,y> = \int_Ul(x)\cup y.
$$
This gives the result.
\end{pf}
\subsection{Hodge structure of projective bundles}

Deligne's construction  \cite{deligne-hodge} gives a stronger result
than is actually stated, namely that for a smooth quasi-projective
variety $X$, $H^i(X)$ take values in the subcategory of polarizable mixed
Hodge structures \cite{beilinson}. One pleasant feature of
this subcategory is the following generalization of Poincar\'e
reducibility:

\begin{lem} The category of polarizable rational pure  Hodge
  structures  is semisimple.
\end{lem}

\begin{pf} This follows from \cite{beilinson}.
\end{pf}

\begin{cor} The category of polarizable rational pure Hodge structures
 satisfies cancellation, i.e. if $A\oplus B \cong
A\oplus C$ then $B \cong C$.
\end{cor}

\begin{prop}\label{prop:hbundle}
Let $p:M\to N$ be a $(d-1)$-fold fiber product
of $\bbP^{n-1}$-bundles (which need not
be locally trivial
in the Zariski topology) over a smooth simply connected 
quasiprojective variety $N$. Then
\begin{enumerate}
\item[(a)] $H^1(M,\bbZ) = 0$\,;

\item[(b)] $p^* : H^3(N, \bbZ) \to H^3(M,\bbZ)$ is surjective, and its kernel
is finite. In particular $p^*$ induces an isomorphism
$H^3(N,\,\bbZ)_{free} \cong
H^3(M,\,\bbZ)_{free}$ (where $A_{free} = A/A_{tors}$)

\item[(c)] $ H^i(M,\bbQ) = \bigoplus_{j\ge 0}  H^{i-2j}(N, \bbQ)\otimes
H^{2j}((\bbP^{n-1})^{d-1},\bbQ)$ as mixed Hodge structures.

\end{enumerate}

\end{prop}

\begin{pf}
We will use the Leray spectral sequence with
coefficients in $\bbZ$ and $\bbQ$.
As $N$ is simply connected, $R^ip_*\bbZ$ is the constant
sheaf
$H^i((\bbP^{n-1})^{d-1},\bbZ)_N$. 
The first statement is an immediate consequence
of the Leray spectral sequence since $H^1(N,\, p_*\bbZ) = 0$.

Note that aside from $H^3(N,\, p_*\bbZ)=
H^3(N,\,\bbZ)$, all the other $E_2$ terms that contribute to
$H^3(M,\bbZ)$  vanish. This implies that $H^3(M,\, \bbZ)$ is 
a quotient of $H^3(N,\,\bbZ)$. But it is an isomorphism
after tensoring with $\Q$ by \cite{deligne68}. Since
the quotient $H^3(N,\,\bbZ) \twr H^3(M,\,\bbZ)$ arising from the analysis
of the spectral sequence is precisely $p^*$, therefore the kernel
is finite. This implies that the induced map
$H^3(N,\,\bbZ)_{free} \to 
H^3(M,\,\bbZ)_{free}$ is surjective with trivial kernel.

Suppose that $p: M\to N$ is a $\bbP^{n-1}$ bundle, that
is $d=2$.
Let  $D$ be a ``nonvertical'' divisor class on $M$ i.e. a class
which restricts nontrivially to each fiber of $p$ (for instance
$c_1(\omega_{M/N})$).
 Then the subspace generated by $D^i$  can
be identified with a copy of $\bbQ(-i)$ in $H^{2i}(M)$. Thus
one gets a morphism of mixed Hodge structures,
$$\bigoplus_{j < n} H^{i-2j}(N)(-j) \to H^i(M)$$
given by summing the maps $\alpha \mapsto p^*\alpha\cup D^j$. 
Note that
the restrictions of $D^j$ generate the cohomology of each
fiber. Therefore by the Leray-Hirch theorem  the above
maps  are isomorphisms, and this proves the last statement
in this case. In general, $M$ is a fiber product of $\bbP^{n-1}$
bundles. Let $D_k$ be the pullback of a nonvertical divisor class
from the $k$th factor. Arguing as before, we get isomorphisms
$$\bigoplus_{j_1+\ldots j_{d-1} = j,\, j_i<n}  H^{i-2j}(N)(-j)
 \to H^i(M)$$
obtained by summing the maps 
$$\alpha \mapsto 
p^*\alpha\cup D_1^{j_1}\cup\ldots D_{d-1}^{j_{d-1}}$$

\end{pf}

By combining this with the previous corollary, and using the fact that sub-mixed
Hodge structures of pure Hodge structures are pure, we obtain:

\begin{cor}\label{cor:hbundle}
Let $p_i:M_i\to N_i$ $i=1,2$ be two  $(\bbP^{n-1})^{d-1}$-bundles
satisfying the hypotheses on $p$ in Proposition\,\ref{prop:hbundle}. If $\phi:M_1\to M_2$ is a morphism
of varieties inducing isomorphisms (of rational mixed Hodge structures)
$H^i(M_2,\bbQ)\to H^i(M_1,\bbQ)$ for $i \le k$, and if these mixed Hodge
structures are pure of weight $i$ then for $i \le k$ the Hodge structures
$H^i(N_2,\bbQ)$ and $H^i(N_1,\bbQ)$ are  pure of weight $i$ and are
{\em noncanonically} isomorphic as rational pure Hodge structures of weight $i$.
\end{cor}

\section{Hodge structure on degree one moduli}\label{s:hodgedeg1}

Let $\cu_X(n,1)$ be the moduli space of semi-stable rank $n$
vector bundles of degree $1$ on $X$. There is a smooth morphism
$\det:\cu_X(n,1)\to Pic^1(X)$ given by the determinant, and the
fiber over $L'$ is just $\cs\cu_X(n,L')$.

Atiyah and Bott \cite[9.11]{at-bot} gave generators for
the cohomology ring $H^*(\cu_X(n,1),\bbQ)$ in terms of certain
universal characteristic classes. As  the restriction
map on rational cohomology from $\cu_X(n,1)$ to $\cs\cu_X(n,L')$
is surjective [loc. cit. 9.7], these give generators
for the cohomology of the $\cs\cu_X(n,L')$.
We will describe these generators in a form which is convenient
for us.
Let $E$ be a Poincar\'e vector bundle. Fix a line bundle $L'\in
Pic^1(X)$, let $p_i$ denote the projections
on $X\times \cs\cu_X(n,L')$ and let $c_r$ be the $r$-th Chern class
of $E$ restricted to $X\times{\cs\cu(n,L')}$.
Let $\gamma_{r,i}$ denote the  composition of the
following  morphisms of Hodge structures:
\begin{align*}
H^i(X,\bbQ)(-r+1) & \stackrel{p_1^*}{\lr} H^i(X\times \cs\cu_X(n,L'),\bbQ)(-r+1) \\
H^i(X\times \cs\cu_X(n,L'),\bbQ)(-r+1) & \stackrel{c_r\cup}{\lr}
H^{i+2r}(X\times \cs\cu_X(n,L'),\bbQ)(1) \\
H^{i+2r}(X\times \cs\cu_X(n,L'),\bbQ)(1)
& \stackrel{p_{2*}}{\lr} H^{i+2r-2}(\cs\cu_X(n,L'), \bbQ)
\end{align*}

\begin{thm}{\em (Atiyah-Bott \cite{at-bot})} The ring
  $H^*(\cs\cu_X(n,L'),\bbQ)$
is generated by the images of the maps $\gamma_{r,i}$ for $2\le r\le
n$ and $0\le i\le 2$.
\end{thm}

\begin{cor}\label{cor:indepHodge} For each $i$, any simple summand of  
the Hodge structure
$$
H^i(\cs\cu_X(n,L'),\bbQ)
$$
is, up to Tate twisting, a direct summand of a tensor power of 
$H^1(X)$. 
The Hodge structure on $H^i(\cs\cu_X(n,L'),\bbQ)$ 
is independent of $L'\in Pic^1(X)$
\end{cor}

\begin{pf}
Let $P= Pic^1(X)$.
Let $H$ be the direct sum of all the tensor products
$$H^{i_1}(X,\bbQ)(-r_1+1)\otimes H^{i_2}(X,\bbQ)(-r_2+1)\ldots$$
indexed by  the finite sequences $((i_1,r_1),(i_2,r_2),\ldots)$
with
 $$(i_1+2r_1-2)+(i_2+2r_2-2)+\ldots = i.$$
Then, by the theorem, there is a surjection
$H\to H^i(\cs\cu_X(n,L'),\bbQ)$ given by the product
of $\gamma$'s. This implies the first statement.
The above map extends to a surjection of variations of
Hodge structures $H_P\to R^i\det_*\bbQ$, where the
$H_P$ denotes the constant variation with fiber $H$.
 The second statement now
follows from the theorem of the fixed part \cite[4.1.2]{deligne-hodge}.

\end{pf}

Since the relations among the above generators have recently been
determined by Jeffrey and Kirwan \cite{jef-kir}, it is possible
to make a complete determination of these Hodge structures (over
$\Q$). It is not clear whether the maps $\gamma_{r,i}$
are surjective for integral cohomology, however one does have:

\begin{thm}\label{thm:h1Xtoh3}{\em(Narasimhan-Ramanan \cite{N-R})}
The map
$$
\gamma_{2,1}:H^1(X,\bbZ)(-1)\to H^3(\cs\cu_X(n,L'),\bbZ)
$$
is an isomorphism.
\end{thm}

This is not stated as such, but this is implicit in their proof
of their third theorem.

\section{Main theorems}\label{s:mainthms}

\subsection{Natural polarizations}
We give a proof of the following ``folklore'' theorem: The
theta divisor and its multiples are the only natural polarizations
on the Jacobian. To simplify the discussion, we
work with polarizations in the Hodge theoretic sense.
Let $\pi:\cx\to T$ be a family of
genus $g$ curves over an irreducible base variety. Let $\theta$ be the
standard  polarization on $R^1\pi_*\bbZ$ corresponding
to cup product, and let $\theta'$ be
some other polarization.

\begin{lem}\label{lem:natlpolar}
If the canonical map to moduli space $T\to M_g$ is
dominant, then there exist a positive integer $m$ such that
$\theta' = m\theta$.
\end{lem}

\begin{pf} Let $t_0\in T$ be a base point.
$\theta_{t_0}$ can be viewed as a primitive vector in
$V=\wedge^2 H^1(X_{t_0},\bbZ)^{\pi_1(T,t_0)}$.
Therefore it is enough to prove that $V\otimes \bbR$ is
one dimensional.  After replacing $T$ by a Zariski  open subset,
we can assume that
the image $T'\subset M_g$ is disjoint from the locus of curves
with automorphisms, and that $T\to T'$ is smooth. This
guarantees (by Teichmuller theory) that $\pi_1(T')$ surjects onto
the mapping class group which surjects onto  $Sp_{2g}(\bbZ)$,
and furthermore that
the image of $\pi_1(T)$ has finite index in $\pi_1(T')$.
Therefore $\pi_1(T)$ has Zariski dense image in $Sp_{2g}(\bbR)$.
Consequently $V\otimes\bbR = \wedge^2H^1(X_{t_0},\bbR)^{Sp_{2g}(\bbR)}$
which is well known to be spanned by the standard symplectic form.
\end{pf}

\subsection{First main theorem}

\begin{thm} \label{thm:re-main1} Let 
$\i(n,g) = 2(n-1)g-(n-1)(n^2+3n+1)-7$.
Let $X$ be a curve of genus $g\ge 2$. If $n\ge 4$ and $i < \i(n,g)$ are 
integers, then for any pair of line bundles $L, L'$ on $X$,the
mixed Hodge structures  $H^i(\cs\cu_X^s(n,L),\Q)$
and  $H^i(\cs\cu_X^s(n,L'),\Q)$ are
{\em (noncanonically)} isomorphic and  are both  pure of weight $i$.
\end{thm}

\begin{pf} 
There is no loss of generality in assuming that $\deg\, L' = 1$.
In fact, by \ref{cor:indepHodge}, we can assume that $L'$ is a
specific line bundle, namely $L' = L(-D)$, where $D= x^1+\ldots
x^{d-1}$.
Then reverting to the our earlier notation, we have $\cs_1 =
\cs\cu_X(n,L')$ and $\cs =\cs\cu_X(n,L)$.
Consider the diagram
$$\cs_1\stackrel{\pi}{\leftarrow}
 \bbP \supset U \stackrel{f}{\to} U'\subset \cs^s.$$
Then corollaries \ref{cor:keyest1}, \ref{cor:keyest2}  and
Lemma\,\ref{lem:codim} implies
$$H^i(\bbP,\bbQ) \cong H^i(U, \bbQ)$$
$$H^i(U',\bbQ) \cong H^i(\cs^s, \bbQ)$$
for all $i <\i(n,g) $.
The theorem follows from Corollary\,\ref{cor:hbundle} (see also 
Remark\,\ref{rem:fib-prod}).

\end{pf}

\subsection{Second main theorem}

\begin{thm} Let $X$ be a curve of genus $g\ge 2$, $n\ge 2$ an integer
and $L$ a  line bundle  on $X$.  Let $\cs^s=\cs\cu_X^s(n,\,L)$.
\begin{enumerate}
\item[(a)] If $g > {3\over n-1} + {n^2 + 3n + 3\over 2}$
and $n\ge 4$, then
$H^3(\cs^s,\bbZ)$ is a pure Hodge structure of type
$\{(1,2),\,(2,1)\}$, and it carries a
natural polarization making the intermediate Jacobian
$$J^2(\cs^s) = \frac{H^3(\cs^s,\,\bbC)}
{F^2+H^3(\cs^s,\,\bbZ)} $$
into a principally polarized abelian variety. There is an
isomorphism of principally polarized abelian varieties $J(X)\simeq
J^2(\cs^s)$.
\item[(b)] If $\deg L$ is a multiple of $n$, then the conclusions of {\em (a)}
are true for $g\ge 3$, $n\ge 2$ except the case $g=3,n=2$.
\end{enumerate}
\end{thm}

\begin{pf} We concentrate on part (a). The proof of part (b)
is identical if we take the Hecke correspondence of \S\,\ref{s:hecke0},
and the codimension estimates of $\bbP\setminus\bbP^s$ and
$\cs\setminus\cs^s$ given in that section (see \ref{ss:codim}).
Consider the diagram
$$\cs_1\stackrel{\pi}{\leftarrow}
 \bbP \supset U \stackrel{f}{\to} U'\subset \cs^s.$$
once again. Then corollaries \ref{cor:keyest1}, \ref{cor:keyest2}  and
Lemma \ref{lem:codim} imply
\begin{align*}
H^3(\bbP,\bbZ)& \cong H^3(U,\,\bbZ), \\
H^3(U',\bbZ) & \cong H^3(\cs^s,\,\bbZ).
\end{align*}
 $\cs_1$ is unirational 
(see \cite[pp.\,52--53,\,VI.B]{css-drez}),
 hence so is $\bbP$. Therefore  these varieties
are simply connected \cite{pi1}. Since
$\codim(\bbP\sm U) > 1$ it follows that $U$ is simply connected (purity
of branch locus). The homotopy exact sequence tells us that $U'$ is simply
connected.
Proposition \ref{prop:hbundle} applied to $\pi$ and $f$ implies that
on the third cohomology, $\pi^*$ and $f^*$ are surjective with finite
kernels. Since $\pi$ is locally trivial in the Zariski topology, $\pi^*$
is in fact an isomorphism. Combining these facts with the isomorphisms above
produces a map of Hodge structures
$$H^3(\cs^s,\,\bbZ)\to H^3(\cs_1,\,\bbZ)$$ which is surjective with
finite kernel. As an immediate consequence we have,
$$H^3(\cs_1,\bbZ)_{free}\cong H^3(\cs^s, \bbZ)_{free} .$$
This, together with with Theorem\,\ref{thm:h1Xtoh3}, yields
an isomorphism of Hodge structures:
\begin{equation}\label{eq:h1Xtoh3}
H^1(X,\bbZ)(-1) \cong H^3(\cs^s, \bbZ)_{free} 
\end{equation}
which yields an isomorphism of tori $J(X)\cong J^2(\cs^s)$.

The next step is to construct a polarization on $H^3(\cs^s, \bbZ)$.
One knows from the results of Drezet and Narasimhan \cite{drez-nar}, that
$\pic(\cs^s)=\bbZ$ (see p.\,89, 7.12 (especially the proof)
of {\it loc.cit.}). Moreover,
$\pic(\cs) \to \pic(\cs^s)$ is an isomorphism. Let $\xi'$ be the ample
generator of $\pic(\cs^s)$. It is easy to see that there exists a positive
integer $r$, independent of $(X,\,L)$ (with genus $X=g$), such that
$\xi = {\xi'}^r$ is very ample on $\cs$ (we are not distinguishing
between line bundles on $\cs^s$ and their (unique) extensions to
$\cs$). Embed $\cs$ in a suitable projective space via $\xi$. Let
$e=\codim(\cs\sm \,\cs^s)$. Let $M$ be the intersection of $k = \dim\cs -e +1$
hyperplanes (in general position) with $\cs^s$. Then $M$ is smooth, projective
and contained in $\cs^s$. Let $p = \dim{\cs}$ and $H^*_c$ --- cohomology with
compact support. We then have a map
$$
l\col H^3(\cs^s) \lr H^{2p-3}_c(\cs^s)
$$
defined by
$$
x \mapsto x\cup c_1(\xi)^{p-k-3}\cup [M].
$$
If $M'$ is another $k$-fold intersection of general hyperplanes, then
$[M'] = [M]$. Hence $l$ depends only on $\xi$. According to
Proposition\,\ref{prop:polarization},
the pairing on $H^3(\cs^s,\,\bbC)$ given by
$$
<x,\,y> = \int_{\cs^s}l(x)\cup\,y
$$
gives a polarization on the Hodge structure of $H^3(\cs^s)$.
Pulling this back via the isomorphism \eqref{eq:h1Xtoh3} gives
a second polarization on $H^1(X)$. If we can show that $<,>$
varies over the whole $M_g$, then we can appeal to
lemma \ref{lem:natlpolar} to show that there exists a
positive integer $m$ (independent of $X$)
such that $\frac{1}{m}<,>$ coincides with the standard principal
polarization of $H^1(X)$, and this will complete the proof.

Let $T$ be the moduli space parameterizing tuples
$$(X,x^1,\ldots x^{d-1}, L,\lambda)$$
consisting of a genus $g$ curve, $d-1$ distinct points, a degree $d$
line bundle and a level $3$ structure. $T$ is an irreducible
variety which surjects onto $M_g$. The inclusion of
the level structure
guarantees that $T$ is fine, and therefore there is a universal
curve $\cx\to T$ together with $d-1$ sections and so on.
All of the constructions given so far can be carried out relative
to $T$, in other words we can construct a diagram of $T$-schemes
$$
\begin{array}{ccccc}
S_1&\leftarrow  & P          & \dashrightarrow& S \\
   &\searrow p_1& \downarrow & p\swarrow      & \\
   &            &  T         &                 &
\end{array}
$$
whose fibers are the Hecke correspondences. Furthermore there
is  polarization on $R^3p_*\bbZ\cong R^1p_{1*}\bbZ(-1)$,
which restricts to the one constructed in the
previous paragraph. We now in a position to apply lemma
\ref{lem:natlpolar}
to conclude the proof.

\end{pf}

\begin{ack} We wish to thank Prof.\,M.\,S. Narasimhan, Prof.\,C.\,S.
Seshadri, and Prof.\,S. Ramanan for their encouragement and their help. 
Thanks to V. Balaji, L. Lempert, N. Raghavendra and P.\,A.Vishwanath for 
helpful discussions.  Balaji made us aware of the problem, and generously 
discussed his proof (in \cite{balaji}) of the Torelli theorem for Seshadri's 
desingularization of $\cs\cu_X(2,\,\co_X)$. We wish to thank the anonymous
referee for the care with which he/she read the paper, and for the helpful
comments.  The second author gratefully acknowledges
the four wonderful years he spent at the SPIC Science Foundation, Madras.
\end{ack}


\bibliographystyle{plain}
\end{document}